# Posterior consistency of Dirichlet mixtures of beta densities in estimating positive false discovery rates

**Subhashis Ghosal**[*1], **Anindya Roy**[†2] and **Yongqiang Tang**[3]

*North Carolina State University, University of Maryland – Baltimore County and SUNY Downstate Medical Center*

**Abstract:** In recent years, multiple hypothesis testing has come to the forefront of statistical research, ostensibly in relation to applications in genomics and some other emerging fields. The false discovery rate (FDR) and its variants provide very important notions of errors in this context comparable to the role of error probabilities in classical testing problems. Accurate estimation of positive FDR (pFDR), a variant of the FDR, is essential in assessing and controlling this measure. In a recent paper, the authors proposed a model-based nonparametric Bayesian method of estimation of the pFDR function. In particular, the density of p-values was modeled as a mixture of decreasing beta densities and an appropriate Dirichlet process was considered as a prior on the mixing measure. The resulting procedure was shown to work well in simulations. In this paper, we provide some theoretical results in support of the beta mixture model for the density of p-values, and show that, under appropriate conditions, the resulting posterior is consistent as the number of hypotheses grows to infinity.

## 1. Introduction

Consider the problem of testing $m$ null hypotheses $H_{0,1}, \ldots, H_{0,m}$ simultaneously, where $m$ is a large number. This type of multiple hypothesis testing problem has received a lot of attention in recent years, primarily due to advanced data collection techniques in genomics, microarray analysis, proteomics, fMRI and some other fields. The analog of type I error probability in multiple testing problems is given by the family-wise error rate, which is defined as the probability of making at least one false rejection. Such a measure is too stringent when $m$ is even moderately large and will block many genuine discoveries (i.e., rejection of a false null hypothesis). In a pioneering paper, Benjamini and Hochberg [2] introduced the concept of the false discovery rate (FDR), the expected value of the ratio of the number of false rejections to the total number of rejections, and described a procedure to control it. Mathematically, the FDR at a nominal level $\gamma$ is given by $\mathrm{E}(V/\max(R,1)) = \mathrm{E}(V/R|R > 0)\mathrm{P}(R > 0)$, where $R = R(\gamma)$ stands for the number of hypotheses rejected at nominal level $\gamma$ and $V = V(\gamma)$ is the number of false rejections among these. Storey [11, 12] argued that the positive false discovery rate (pFDR) (at nominal level $\gamma$) defined as $\mathrm{E}(V/R|R > 0)$, is a more relevant measure to control. Storey's approach consists of estimating the pFDR function at each $\gamma$ and choosing a $\gamma$ so that the estimated pFDR function is within a given limit, $\alpha$. Storey [11, 12] showed that under a certain natural setup, the resulting procedure controls pFDR by $\alpha$. Some other related measures have also been considered in the literature; see Benjamini and Hochberg [2], Efron and Tibshirani [3], Tsai et al. [14] and Sarkar [10].

---
*Supported in part by NSF Grant DMS-03-49111.
†Supported in part by NIH Grant 1R01GM075298-01.
[1]Department of Statistics, North Carolina State University, 2501 Founders Drive, Raleigh, NC 27695, USA, e-mail: ghosal@stat.ncsu.edu
[2]Department of Mathematics and Statistics, University of Maryland Baltimore County, 1000 Hilltop Circle, Baltimore, MD 21250, USA, e-mail: anindya@math.umbc.edu
[3]SUNY Downstate Medical Center, 450 Clarkson Avenue, Brooklyn, NY 11203, USA, e-mail: yongqiang_tang@yahoo.com

In order to estimate the pFDR function, Storey [11] considered a mixture model setup, where each null hypothesis has a fixed probability, $\pi$, of being true. Thus, the number of true null hypotheses, $m_0$, is taken to be a random variable distributed as binomial $(m, \pi)$. If the null hypothesis is true, then it is assumed that the p-value associated with the corresponding test statistic is uniformly distributed. The p-value when the alternative is true and has a fixed value $\theta$, follows a distribution $H = H(\cdot|\theta)$. It is somewhat unnatural to assume that the alternative value remains fixed when the hypotheses themselves are appearing randomly. A more natural assumption would be to assume that, given that null hypothesis is false, the alternative is chosen randomly according to a distribution $\mu$. Then, marginally, the conclusion that the p-value under the alternative is distributed as $H$ remains unaffected, where now $H$ stands for the mixture $\int H(\cdot|\theta)d\mu(\theta)$. Under this setup, Storey [11] showed that the pFDR at nominal level $\gamma$ is given by the expression $\pi\gamma/[\pi\gamma + (1 - \pi)H(\gamma)]$. To estimate the pFDR, it then suffices to estimate $\pi$, since the denominator can be estimated essentially by the empirical proportion $R/m$. Actually, Storey [11] used a slightly different estimator to take into account the problem of zeros in finite samples. Estimation of $\pi$ is more delicate. Storey [11] assumed that for some appropriate threshold value $\lambda$, all p-values over $\lambda$ are associated with true null hypotheses. Equating the observed proportion of rejected hypotheses with the expectation $\lambda(1 - \pi)$, and choosing $\lambda$ appropriately, an estimate of $\pi$, and hence that of pFDR, can be obtained.

Although Storey [11] did not make any explicit assumption about $H$, implicitly it was assumed that $H$ is concentrated near zero. It is this assumption that leads to the conclusion that almost every p-value over level $\lambda$ must arise from null hypotheses. While this is reasonable, it introduces some bias in the analysis because, although relatively rare, p-values bigger than $\lambda$ can occur under alternatives as well.

The density of p-values under alternatives usually has more features than is assumed above. These important features may be exploited to construct a more refined estimator of pFDR. For instance, the density of p-values under an alternative value is often decreasing, dropping from an infinite height at 0 to a very low or no height at 1, and the derivative of the density approaches zero near the point 1. These densities resemble beta $(a, b)$ densities $\mathrm{be}(x; a, b) = (1/B(a, b))x^{a-1}(1-x)^{b-1}$ with $a < 1$ and $b \geq 1$, or their mixtures, where $B(a, b) = \Gamma(a)\Gamma(b)/\Gamma(a+b)$ is the beta function. A reasonable model may be proposed for this type of densities, and based on the model it may be possible to estimate the pFDR function more accurately. Tang et al. [13] modeled the p-value density under the alternative as a mixture of beta densities and thereby incorporated some of the salient features of the p-value density directly into the model. They followed a Bayesian approach by putting a Dirichlet process prior on the mixing distribution of the beta parameters. The resulting posterior is amenable to Markov chain Monte-Carlo methods of com-putation. Tang et al. [13] showed by simulation that the resulting procedure gives more stable and accurate estimates of the pFDR function.

In this paper, we theoretically study the appropriateness of the model assumptions made in Tang et al. [13] and investigate the support of the Dirichlet mixture of beta prior. Our results provide important theoretical justification for the setup assumed in Tang et al. [13]. Under certain conditions, we show that the posterior distribution of the pFDR function is consistent as the number of hypotheses tends to infinity.

## 2. Mixture model framework

### *2.1. Basic setup*

Suppose we have observed the values of the test statistics for testing $m$ null hypotheses $H_{0,i}$, $i = 1, \ldots, m$, against appropriate alternatives. Let $X_1, \ldots, X_m$ stand for the p-values for the respective $m$ tests. We assume that the tests are based on independent data, so that $X_1, \ldots, X_m$ are independent. We also assume that there is a random mechanism which independently determines whether $H_{0,i}$'s are true or false, respectively with probability $\pi$ and $1 - \pi$. Let $H_i = I(H_{0,i}$ is true), be the indicator that the $i$th null hypothesis is true. Of course, $H_i$'s are unobserved.

The distribution of $X_i$ under $H_{0,i}$ can be assumed to be the uniform distribution on $[0, 1]$. This happens whenever the test statistic is a continuous random variable and the null hypothesis is simple, or in situations like the t-test or F-test, where the null hypothesis has been reduced to a simple one by considerations of similarity or invariance. Under more general situations, the property can still be expected to be approximately true if, for instance, a conditional predictive p-value or a partial predictive p-value (Bayarri and Berger [1]) is used; see Robins et al. [9] for details. If the null and alternative hypotheses are one-sided and the underlying distribution has the monotone likelihood ratio (MLR) property, then the power function is increasing in the parameter, and, as a result, the null distribution of the p-value is stochastically larger than the uniform. Many estimation procedures remain valid in a conservative sense when the actual null distribution is replaced by the uniform. It is easy to show that Storey's estimators have this property. The Bayesian estimator of Tang et al. [13] also enjoys the same property – see Tang et al. [13] for discussion. Henceforth we shall assume that the null distribution of p-values is $U[0, 1]$.

Let $f(x)$ stand for the density of the p-value under an alternative distribution. The following result shows that under a natural condition, $f(x)$ is decreasing.

**Proposition 1.** *Suppose that the p-value is computed using a statistic, $T$, whose density, $g_\theta$, has the MLR property. Then the p-value density $f(x)$ is decreasing.*

*Proof.* Let $\theta_0$ stand for the value of the parameter under the null hypothesis and $\theta_1$ stand for the value under the alternative. Let $T_{\mathrm{obs}}$ stand for the observed value of $T$. Denote the cumulative distribution function (c.d.f.) of $g_\theta$ by $G_\theta$. Then the distribution function of the p-value under $\theta_1$ is

$$F_{\theta_1}(x) = P_{\theta_1}(P_{\theta_0}(T_n > T_{\mathrm{obs}}) \le x) = 1 - G_{\theta_1}(G_{\theta_0}^{-1}(1 - x)).$$

Hence the p-value density is given by

$$(2.1) \qquad f_{\theta_1}(x) = \frac{g_{\theta_1}(G_{\theta_0}^{-1}(1 - x))}{g_{\theta_0}(G_{\theta_0}^{-1}(1 - x))} = \frac{g_{\theta_1}(z)}{g_{\theta_0}(z)},$$

where $z = G_{\theta_0}^{-1}(1-x)$. By the MLR property, the expression in (2.1) is increasing in $z$, equivalently decreasing in $x$. $\qquad\square$

For standard two-sided tests like the $z$-test or $t$-test, the density of the p-value under the alternative is also decreasing. Under certain assumptions which are satisfied generally, the following result shows a two-sided analog of the previous proposition.

**Proposition 2.** *Suppose that the p-value is computed using a statistic $T$ whose density $g_\theta$ is symmetric under the null hypothesis $H_0 : \theta = \theta_0$. Further suppose that for the symmetrized density $\tilde{g}_\theta(z) = (g_\theta(z) + g_\theta(-z))/2$, the ratio $\tilde{g}_\theta(z)/g_{\theta_0}(z)$ is increasing in $z$. Then the p-value density $h(x)$ is decreasing.*

*Proof.* With notations as in the last proof, the distribution function of the p-value under $\theta_1$ is

$$F_{\theta_1}(x) = P_{\theta_1}(2P_{\theta_0}(T_n > |T_{\text{obs}}|) \le x)$$
$$= 1 - G_{\theta_1}(G_{\theta_0}^{-1}(1-x/2)) + G_{\theta_1}(-G_{\theta_0}^{-1}(1-x/2)).$$

The p-value density can be seen to be given by

$$(2.2) \quad f(x) = f_{\theta_1}(x) = \frac{g_{\theta_1}(G_{\theta_0}^{-1}(1-x/2))}{2g_{\theta_0}(G_{\theta_0}^{-1}(1-x/2))} + \frac{g_{\theta_1}(-G_{\theta_0}^{-1}(1-x/2))}{2g_{\theta_0}(-G_{\theta_0}^{-1}(1-x/2))} = \frac{\tilde{g}_{\theta_1}(z)}{g_{\theta_0}(z)},$$

which is decreasing in $x$ by the given assumption. $\qquad\square$

The p-value density for a one-sided hypothesis generally decays to zero as $x$ tends to 1. Let $L$ stand for the lower limit of the value of the test statistic, which is often $-\infty$. Assume that as $z \to L$, we have that $g_{\theta_1}(z)/g_{\theta_0}(z) \to 0$. Then, clearly it follows from (2.1) that $f(x) \to 0$ as $x \to 1$, since $z = G_{\theta_0}^{-1}(1-x) \to L$ as $x \to 1$.

For a two-sided hypothesis, $\tilde{g}_{\theta_0}(z)/g_{\theta_0}(z)$ will not generally go to 0 as $z \to L$, and hence the minimum value of the p-value density will be a (small) positive number. For instance, for the two-sided normal location model, the minimum value is $e^{-n\theta^2/2}$, where $n$ is the sample size on which the test is based.

## 2.2. Identifiability and continuity properties

If a c.d.f. $F$ on [0,1] can be written as $F(x) = \pi x + (1-\pi)H(x)$, where $H(\cdot)$ is another c.d.f. on [0,1], then the representation is generally not unique, so that $\pi$ and $H$ are not separately identifiable. The components $\pi$ and $H$ can be identified by imposing the additional condition that $H$ cannot be represented as a mixture with another uniform component, which, for the case when $H$ has a continuous density $h$, translates into $h(1) = 0$. Define the map $\pi(F)$ from the space of continuous c.d.f. on [0,1] to [0,1] as the maximum possible value of $\pi$ in the mixture representation $F(x) = \pi x + (1-\pi)H(x)$. As in all mixture problems, $H$ is not defined when $\pi(F)$ is one, that is, $F$ is the uniform distribution on [0,1]. When $F$ physically stands for the p-value distribution, $\pi(F)$ is an upper bound for the proportion of null hypothesis and therefore $\pi(F)\gamma/F(\gamma)$ is an upper bound for the actual pFDR. Thus this choice of $\pi$ is appropriate in a conservative sense in that in order to control pFDR, it suffices to control the auxiliary quantity $\overline{\text{pFDR}}(F;\gamma)$ defined as $\pi(F)\gamma/F(\gamma)$.

Let $\mathcal{F}$ stand for all $F$ representable as $F(x) = \pi x + (1-\pi)H(x)$ for $\pi \in [0,1]$. The following proposition shows an important upper-semicontinuity property of the map $\pi(F)$. Let $\to_w$ stand for weak convergence of probability distributions.

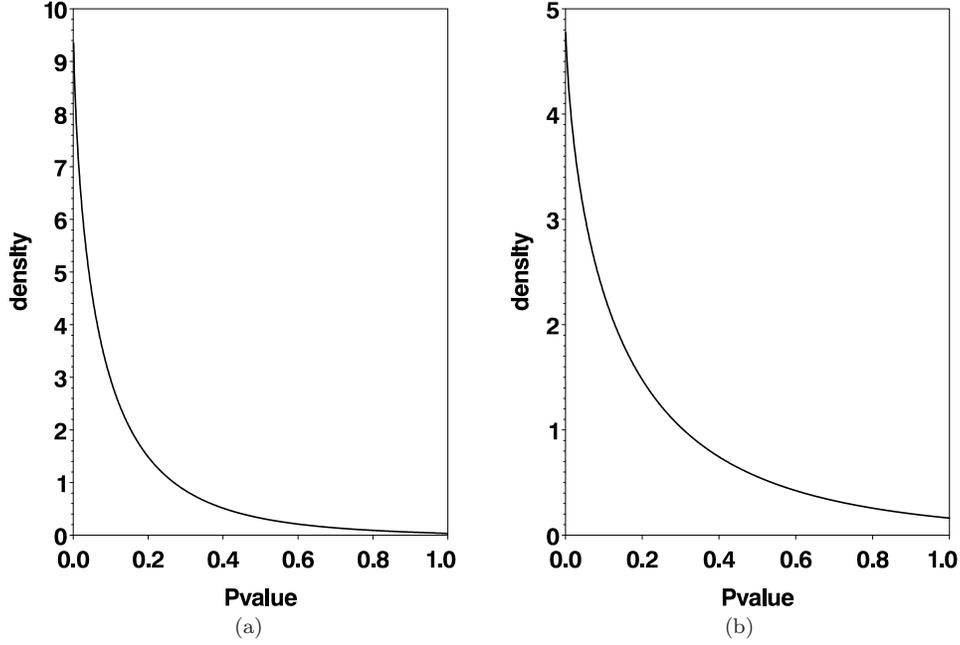

FIG 1. *Plots of p-value density for t-test with 3 d.f. (a) p-value density of one-sided t-test. (b) p-value density of two-sided t-test.*

**Proposition 3.** *The class $\mathcal{F}$ is weakly closed and the map $F \mapsto \pi(F)$ on $\mathcal{F}$ is upper-semicontinuous, that is,*

$$\text{if } F_n \to_w F, \text{ then } \limsup_{n\to\infty} \pi(F_n) \leq \pi(F).$$

*Further, for any $\gamma$, $\limsup_{n\to\infty} \overline{\mathrm{pFDR}}(F_n; \gamma) \leq \overline{\mathrm{pFDR}}(F; \gamma)$.*

*Proof.* Let $F_n \in \mathcal{F}$ and $F_n \to_w F$. Because $\pi_n = \pi(F_n)$ is a bounded sequence and $H_n$ in the representation $F_n(x) = \pi_n x + (1 - \pi_n) H_n(x)$ is tight, we may assume that both are convergent along a subsequence, to $\pi^*$ and $H^*$, respectively. Then $F(x) = \pi^* x + (1 - \pi^*) H^*(x)$, and hence $F \in \mathcal{F}$.

Observe that for any $F \in \mathcal{F}$, $\bar{F}(\lambda) \geq \pi(F)(1 - \lambda)$ for all $0 \leq \lambda \leq 1$ and that $\pi(F) = \inf\{\bar{F}(\lambda)/(1 - \lambda) : 0 < \lambda \leq 1\}$. The infimum is attained because, by our choice, $\pi(F)$ is the largest $\pi$ in the mixture representation.

Now for any fixed $\lambda_0$ which is a continuity point of $F$, we have that

$$\limsup_{n\to\infty} \pi(F_n) = \limsup_{n\to\infty} \inf_{\lambda} \frac{\bar{F}_n(\lambda)}{1 - \lambda} \leq \lim_{n\to\infty} \frac{\bar{F}_n(\lambda_0)}{1 - \lambda_0} = \frac{\bar{F}(\lambda_0)}{1 - \lambda_0}.$$

Since $\lambda_0$ is arbitrary and the set of continuity points of $F$ is dense in [0,1], the first assertion follows.

The last relation clearly follows from the expression for $\overline{\mathrm{pFDR}}$. $\square$

Under additional restrictions, identifiability of the components $\pi$ and $H$ and continuity of $\pi(F)$ may be established. For example, the following class of c.d.f. $F$ allows $\pi$ and $H$ to be identified from $F$. Assume that the p-value distribution $H$ under the alternative belongs to $\mathcal{D}$, the class of c.d.f. on [0,1] which admits a density $h$, with $h(1) = 0$. Let $\mathcal{F}_{\mathcal{D}}$ denote the class of all c.d.f. on [0,1] of the form

$F(x) = \pi x + (1 - \pi) H(x)$ for $\pi \in (0, 1)$ and $H \in \mathcal{D}$. Let $f_{\pi, h} = \pi + (1 - \pi) h$ be the corresponding mixture density.

**Proposition 4.** *If $f_{\pi, h} = f_{\pi^*, h^*}$, then $\pi = \pi^*$ and $h = h^*$.*

*Proof.* $f_{\pi, h} = f_{\pi^*, h^*}$ implies $\pi + (1 - \pi) h(x) = \pi^* + (1 - \pi^*) h^*(x)$ for all $x$. Putting $x = 1$ and using the fact that $h(x) = h^*(x) = 0$, we have $\pi = \pi^*$. This now implies $h = h^*$ or $H = H^*$. ◻

To study consistency, we need to show that $\pi$ and $h$ can be continuously solved from $f$. However, the class $\mathcal{F}_{\mathcal{D}}$ is not weakly closed. We need to impose a restriction on the class of alternative densities so that the tail at 1 remains thin even in the weak limit. Let $\mathcal{B}$ denote a class of c.d.f. on [0,1] that is weakly closed and for all $H \in \mathcal{B}$ we have $\lim_{y \to 0} y^{-1} \bar{H}(1 - y) = 0$. The interval $(1 - y, 1]$ is open in $[0, 1]$. Hence, by the portmanteau theorem, $H_n \to_w H$ implies that $\bar{H}(1 - y) \le \liminf_{n \to \infty} \bar{H}_n(1 - y)$. Thus for the weak limit $H$ of a sequence $H_n \in \mathcal{B}$ to be in $\mathcal{B}$, one needs to be able to interchange the order of the limits with respect to $y$ and $n$. For instance, if $\mathcal{B} = \{H : \bar{H}(1 - x) \le \psi(x) \text{ for all } x < \delta\}$, where $\delta > 0$ is a fixed number and $\psi$ is a fixed function which satisfies $\psi(x) = o(x)$ as $x \to 0$ (like $Cx^{1+\epsilon}$), then the class $\mathcal{B}$ satisfies the requirement. Let $\mathcal{F}_{\mathcal{B}}$ denote the class of c.d.f. on [0,1] representable as $F(x) = \pi x + (1 - \pi) H(x)$ for $\pi \in (0, 1)$ and $H \in \mathcal{B}$. Note that $\mathcal{F}_{\mathcal{B}}$ need not be a subset of $\mathcal{F}_{\mathcal{D}}$ as the c.d.f. in $\mathcal{B}$ need not have a density.

**Proposition 5.** *Identifiability in Proposition 4 holds if $F \in \mathcal{F}_{\mathcal{B}}$.*

*Proof.* If $\pi x + (1 - \pi) H(x) = \pi^* x + (1 - \pi^*) H^*(x)$ for all $x$, then

$$\pi(1 - x) + (1 - \pi) \bar{H}(x) = \pi^*(1 - x) + (1 - \pi^*) \bar{H}^*(x).$$

Dividing both sides by $1 - x$ and letting $x \to 1$, we obtain $\pi = \pi^*$ and hence $H = H^*$. ◻

**Proposition 6.** *The map $(\pi, H) \mapsto F_{\pi, H}$ is a homeomorphism from $(0, 1) \times \mathcal{B}$ to $\mathcal{F}_{\mathcal{B}}$, where $\mathcal{B}$ and $\mathcal{F}_{\mathcal{B}}$ are the weak topology.*

*Proof.* (Forward side) If $\pi_n \to \pi$ and $H_n \to_w H$, then $H_n(x) \to H(x)$ at all continuity points $x$, giving $\pi_n x + (1 - \pi_n) H_n(x) \to \pi x + (1 - \pi) H(x)$.

(Reverse side) Let $F_{\pi_n, H_n} \to_w F_{\pi, H}$. To show that $\pi_n \to \pi$ and $H_n \to_w H$. Fix any subsequence $n'$. It is enough to extract a further subsequence $n''$ along which $\pi_{n''} \to \pi$ and $H_{n''} \to_w H$.

Because $\pi_{n'}$ is bounded and $H_{n'}$ is tight, we can extract a further subsequence $n''$ such that $\pi_{n''} \to \pi^*$ and $H_{n''} \to_w H^*$ for some $\pi^*$ and $H^*$. By the closedness of $\mathcal{B}$ under the weak topology, $H^* \in \mathcal{B}$ (note that $(1 - x, 1]$ is an open subset of $[0, 1]$). By the forward side, $F_{\pi_{n''}, H_{n''}} \to_w F_{\pi^*, H^*}$. Thus $F_{\pi^*, H^*} = F_{\pi, H}$. By identifiability in the class $\mathcal{F}_{\mathcal{B}}$, $\pi^* = \pi$ and $H^* = H$, and hence $\pi_{n''} \to \pi$ and $H_{n''} \to_w H$. This completes the proof. ◻

### 2.3. Mixtures of beta densities

The shape of p-value densities under alternatives has similarities with the beta density be$(x; a, b) = (1/B(a, b)) x^{a-1} (1 - x)^{b-1}$, $0 < x < 1$, for $a < 1$ and $b \ge 1$. Indeed, for the exponential model $\lambda e^{-\lambda z}$, $z > 0$, with parameter $\lambda$ and hypotheses $H_0 : \lambda = \lambda_0$ against $H : \lambda > \lambda_0$, it follows from elementary calculations that the p-value density is exactly beta$(a, 1)$ for some $a < 1$. Mixtures of beta $(a, b)$ with $a < 1$ and $b \geq 1$ make up a considerably large class still preserving the shape of the p-value density, and hence can be considered as a model for p-value densities under the alternative. The following result shows that many similar-shaped densities can be pointwise represented as a mixture of beta $(a, 1)$, a much narrower class.

Recall that a function $\varphi$ on $[0, \infty]$ is called *completely monotone* if it has derivatives $\varphi^{(n)}$ of all orders and $(-1)^n \varphi^{(n)}(z) \geq 0$ for all $z \geq 0$ and $n = 1, 2, \ldots$.

**Proposition 7.** *If a density $h(x)$ on $(0, 1)$ with c.d.f. $H$ can be represented as $h(x) = \int_0^1 a x^{a-1} dG(a)$ for all $0 < x < 1$, then $H(e^{-y})$ is a completely monotone function of $y$ on $[0, \infty)$.*

*Conversely, if $h(x)$ is decreasing and $H(e^{-y})$ is completely monotone, then $h(x) = \int_0^\infty a x^{a-1} dG(a)$ for some probability measure $G$ on $(0, \infty)$ with $\int a^2 dG(a) \leq \int a \, dG(a)$.*

*Proof.* If $h(x)$ is a mixture of be$(a, 1)$, we have that

$$H(x) = \int_0^1 x^a dG(a) = \int_0^1 e^{-a \log x^{-1}} dG(a).$$

Thus $H(x)$ is the Laplace transform of $G$ at the point $\log x^{-1}$. Put $y = \log x^{-1}$ so that $x = e^{-y}$ and $H(e^{-y}) = \int_0^1 e^{-ay} dG(a)$, the Laplace transform of the probability measure $G$. Hence it is completely monotone by Theorem 1 of Section XIII.4 of Feller (1971).

To prove the converse, applying the same theorem and using the fact that $H(1) = 1$, we obtain the representation that $H(e^{-y}) = \int_0^\infty e^{-ay} dG(a)$ for some probability measure $G$ on $(0, \infty)$. Thus $H(x) = \int_0^\infty x^a dG(a)$, and so $h(x) = \int_0^\infty a x^{a-1} dG(a)$. Now, as $h$ is decreasing, $0 \geq h'(x) = \int a(a-1) x^{a-2} dG(a)$. The result now follows by letting $x \to 1$. □

Observe that $\int a^2 dG(a) \leq \int a \, dG(a)$ holds if $G$ is concentrated on $(0, 1]$, but it is not necessary.

**Remark 1.** By a similar argument, if a density $h(x)$ on $(0, 1)$ can be represented as $h(x) = \int_1^\infty b(1 - x)^{b-1} dG(b)$ for all $0 < x < 1$, then the function $\bar{H}(1 - e^{-y})$ is completely monotone as a function of $y$, where $\bar{H}(x) = 1 - H(x)$.

Conversely, if $\bar{H}(1 - e^{-y})$ is completely monotone in $y$ and $h(x)$ is decreasing, then $h(x) = \int_0^\infty b(1 - x)^{b-1} dG(b)$ for some probability measure $G$ on $(0, \infty)$ with $\int b d G(b) \leq \int b^2 dG(b)$.

**Proposition 8.** *Let $\mathcal{H}_1$ stand for the class of decreasing densities $h$ such that $H(e^{-y})$ is completely monotone and $\mathcal{H}_2$ stand for the class of decreasing densities $h$ such that $\bar{H}(1 - e^{-y})$ is completely monotone. A density $h(x)$ on $(0, 1)$ can be represented as a mixture of be$(a, b)$ if $h(x)$ is a convex combination of densities of the form $c h_1(x) h_2(x)$ where $h_1 \in \mathcal{H}_1$ and $h_2 \in \mathcal{H}_2$.*

*Proof.* Clearly it suffices to assume that $h(x) = c h_1(x) h_2(x)$, where $h_1(x) = \int_0^\infty a x^{a-1} dG_1(a)$, $h_2(x) = \int_0^\infty b(1 - x)^{b-1} dG_2(b)$ and

$$c^{-1} = \int_0^\infty \int_0^\infty ab B(a, b) dG_1(a) dG_2(b).$$

Now defining $dG(a, b) = cab B(a, b) dG_1(a) dG_2(b)$, we may write $h(x) = \int be(x; a, b) dG(a, b)$. The total mass of $G$ is given by

$$\int_0^\infty \int_0^\infty cab B(a, b) dG_1(a) dG_2(b) = cc^{-1} = 1,$$

so that $G$ is also a probability measure. This completes the proof. ∎

### *2.4. Dirichlet mixture prior*

Tang et al. [13] proposed a Dirichlet process prior (Ferguson [5]) for the mixing distribution $G$. The parameters of a Dirichlet process $\mathrm{DP}(G_0, \tau)$ are the center measure $G_0 = \mathrm{E}(G)$, and the precision parameter $\tau > 0$. The center measure $G_0$ is the subjective guess about $G$, while $\tau$ controls the concentration of $\mathrm{DP}(G_0, \tau)$ around $G_0$.

The equivalent hierarchical representation in terms of latent variable $(a_i, b_i)$,

$$X_i | a_i, b_i \sim \pi + (1 - \pi)\mathrm{be}(x_i | a_i, b_i),$$
$$(a_i, b_i) | G \overset{\text{i.i.d.}}{\sim} G,$$
$$G \sim \mathrm{DP}(G_0, \tau),$$

is extremely useful in developing the relevant MCMC algorithms for the computation of posterior. Tang et al. [13] used the reparameterization $a = \exp(-|L_a|)$ and $b = \exp(|L_b|)$, and specified $G_0(a, b) = N(L_a | 0, \sigma_a^2) N(L_b | 0, \sigma_b^2)$. Actually, any base measure with full support on $(0, 1) \times (1, \infty)$ will lead to a Dirichlet process with large support.

## 3. Asymptotic properties of posterior

Consider a prior $\Pi$ for $H$ and independently a prior $\mu$ for $\pi$ with full support on $[0, 1]$. Let the true value of $\pi$ and $h$ be, respectively, $\pi_0$ and $h_0$ where $0 < \pi_0 < 1$.

**Theorem 1** (General consistency). *If $h_0$ belongs to the $L_1$-support of $\Pi$ in the sense that $\Pi(\|h - h_0\|_1 < \epsilon) > 0$ for all $\epsilon > 0$, then for every $\epsilon > 0$, $\Pr(\sup\{|F(x) - F_0(x)| : 0 \le x \le 1\} < \epsilon | X_1, \ldots, X_m) \to 1$ a.s.*

*Proof.* For any sequence $F_n$ such that $F_n(x) \to F_0(x)$ for all $x$, continuity of $F_0$ and Polya's theorem imply that $\sup_x |F_n(x) - F_0(x)| \to 0$. Thus given $\epsilon > 0$, we can find a weak neighborhood $\mathcal{W}$ of $F_0$ such that $F \in \mathcal{W}$ implies $\sup_x |F(x) - F_0(x)| < \epsilon$. Thus it suffices to prove that for any weak neighborhood $\mathcal{W}$ of $F_0$,

$$\Pr\{|\pi - \pi_0| < \epsilon, F \in \mathcal{W} | X_1, \ldots, X_m\} \to 1 \text{ a.s. as } m \to \infty.$$

By Schwartz's theorem for weak consistency (see Theorem 4.4.2 of Ghosh and Ramamoorthi [8]), it suffices to show that for every $\epsilon > 0$,

$$(\mu \times \Pi) \left\{ (\pi, h) : \int f_{\pi_0, h_0} \log \frac{f_{\pi_0, h_0}}{f_{\pi, h}} < \epsilon \right\} > 0.$$

Now $f_{\pi, h} \ge \pi$, so $f_{\pi_0, h_0} / f_{\pi, h} \le \pi^{-1} f_{\pi_0, h_0}$, which is integrable, and the integral $\pi^{-1}$ is bounded by a constant when $\pi$ lies in a neighborhood of $\pi_0$. So by Lemma 7 of Ghosal and van der Vaart [7] or Theorem 5 of Wong and Shen [15]

$$\int f_{\pi_0, h_0} \log \frac{f_{\pi_0, h_0}}{f_{\pi, h}} \le A d_H^2(f_{\pi_0, h_0}, f_{\pi, h}) \log_+ \frac{1}{d_H^2(f_{\pi_0, h_0}, f_{\pi, h})},$$

where $d_H$ stands for the Hellinger distance. Also, as $d_H^2(f, g) \le \|f - g\|_1$, it suffices to show that $L_1$-neighborhoods of $f_{\pi_0, h_0}$ gets positive probabilities under $\mu \times \Pi$.

Now,

$$\int_0^1 |[\pi + (1-\pi)h(x)] - [\pi_0 + (1-\pi_0)h_0(x)]|dx$$

$$\leq |\pi - \pi_0| + \int_0^1 |(1-\pi) - (1-\pi_0)|h(x)dx$$

$$+ (1-\pi_0)\int_0^1 |h(x) - h_0(x)|dx$$

$$\leq 2|\pi - \pi_0| + \|h - h_0\|_1.$$

Since $\mu$ gives positive probabilities to neighborhoods of $\pi_0$ and $\Pi$ gives positive probabilities to $L_1$-neighborhoods of $h_0$, the condition of prior positivity holds. □

In view of Proposition 3, the following "upper semi-consistency" (a form of a one-sided consistency) may be concluded.

**Corollary 1.** *Under the conditions of Theorem 1, we have that for any $\epsilon > 0$, $\Pr(\pi < \pi_0 + \epsilon|X_1, \ldots, X_n) \to 1$ a.s. and that the posterior mean $\hat{\pi}_m$ satisfies $\limsup_{m\to\infty} \hat{\pi}_m \leq \pi_0$ a.s.*

Unfortunately, the above corollary has limited significance since typically one would not like to underestimate the true $\pi$ (and the pFDR) while overestimation is less serious. In order to ensure that the convergence takes place, we need to enforce additional restriction on the support of the prior to ensure continuity of $\pi(F)$ with respect to the weak topology on the restricted space.

**Corollary 2.** *Assume that $\Pi$ is supported in $\mathcal{B} \cap \mathcal{D}$ and that $h_0$ belongs to the $L_1$-support of $\Pi$. Then for any $\epsilon > 0$, $\Pr(|\pi - \pi_0| < \epsilon|X_1, \ldots, X_n) \to 1$ a.s. and that $\hat{\pi}_m \to \pi_0$ a.s.*

*Further, for any $0 < \alpha < 1$ and $\epsilon > 0$,*

$$\Pr\left\{\left|\frac{\pi\alpha}{F(\alpha)} - \frac{\pi_0\alpha}{F_0(\alpha)}\right| < \epsilon|X_1, \ldots, X_n\right\} \to 1 \ a.s.$$

*and the above convergence is uniform for $\alpha$ lying in compact subsets of $(0, 1]$.*

*Proof.* The proof of the first assertion follows from Theorem 1 and Proposition 6.

The second assertion follows from the first because $\pi_n \to \pi_0$ and $F_n(\alpha) \to F_0(\alpha)$ implies that $\pi_n\alpha/F_n(\alpha) \to \pi_0\alpha/F_0(\alpha)$, whenever $0 < F_0(\alpha) < 1$, and this holds whenever $0 < \alpha < 1$. In fact, the convergence is uniform over compact subsets of $(0, 1]$, because $F_0(\alpha)$ remains uniformly bounded below there. □

Now we consider a concrete prior obtained from a Dirichlet mixture of betas: Let $h(x) = \int \text{be}(x; a, b)dG(a, b)$, where $G \sim \text{DP}(\tau, G_0)$ and $G_0$ is a probability measure on $(0, 1) \times (1 + \epsilon, \infty)$ with full support. The lower bound $b \geq 1 + \epsilon$ ensures that

$$\bar{H}(1-x) = \int_{1-x}^1 \frac{1}{B(a, b)}y^{a-1}(1-y)^{b-1}dy$$

$$\leq \int_{1-x}^1 b(1-y)^{b-1}dy = x^b \leq x^{1+\epsilon}$$

since $\text{be}(a, b)$ is stochastically dominated by $\text{be}(1, b)$ (by the MLR property of beta distribution) and taking mixtures preserves bounds for the probability of a given set. This ensures that any $H$ in the support of the prior lies in $\mathcal{B}$. This leads to the following consistency result for a Dirichlet mixture of beta prior.

**Theorem 2** (Full $L_1$-support of beta mixture prior). *For any true $h_0 \in \mathcal{B} \cap \mathcal{D}$ lying in the $L_1$-closure of the above beta mixtures, consistency of pFDR holds for the Dirichlet mixture of beta prior if the center measure $G_0$ has support $[0,1] \times [1+\epsilon, \infty)$.*

*Proof.* First let $h_0(x) = h_{Q_0}(x) = \int \mathrm{be}(x; a, b) dQ_0(a, b)$. Given $\epsilon > 0$, find $\eta > 0$ and $M < \infty$ such that $Q_0\{a < \eta$ or $b > M\} < \epsilon$. Let $Q_0^*$ be $Q_0$ restricted and re-normalized to $[\eta, 1] \times [1, M]$. Then by Lemma A.3 of Ghosal and van der Vaart (2001), it follows that $\|h_{Q_0} - h_{Q_0^*}\|_1 < 2\epsilon$. Thus it suffices to assume that $Q_0$ is supported over $[\eta, 1] \times [1+\epsilon, M]$ for some $\eta > 0$ and $M < \infty$. Now if $Q_n$ is a sequence converging weakly to $Q_0$, we may also assume that $Q_n\{a < \eta$ or $b > M\} < \epsilon$ for all $n$ and so that $\|h_{Q_n} - h_{Q_n^*}\|_1 < 2\epsilon$ and $Q_n^*$ converges weakly to $Q_0$. For any $0 < x < 1$, the beta kernel is a bounded continuous function on $[\eta, 1) \times (1, M]$, and hence $h_{Q_n^*}(x) \to h_{Q_0}(x)$. Scheffe's theorem then implies that $\|h_{Q_n^*} - h_{Q_0^*}\|_1 \to 0$. Thus, given any $\epsilon > 0$, if $Q$ lies in a sufficiently small weak neighborhood of $Q_0$, then $\|h_Q - h_{Q_0}\|_1 < \epsilon$. As the center measure $G_0$ has support $[0,1] \times [1+\epsilon, \infty]$, the corresponding Dirichlet process has full weak support. Thus $h_0$ belongs to the $L_1$-support of the prior, and hence consistency holds by Corollary 2.

Now more generally, if $h_0$ can be approximated by beta mixtures in the $L_1$-sense, then also $h_0$ lies in the $L_1$-support as the support is a closed set. Hence consistency is obtained. ∎

**Remark 2.** Proposition 8 gives a sufficient condition for $h_0$ to be in the $L_1$-closure of beta mixtures.

**Remark 3.** By Fubini's theorem, the result continues to hold even if $\tau$ is given a prior and $G_0$ contains hyperparameters.

## 4. Conclusion

A mixture of beta densities $\mathrm{be}(a, b)$ with $a < 1$ and $b > 1$ forms a rich class of densities with shapes like a reflected J. It is shown that, under various natural scenarios, such densities are appropriate for modeling the density of p-values arising from alternative hypotheses. We have also shown that if for any c.d.f. $H$, $H(e^{-y})$ is a completely monotone function of $y$, then the corresponding density $H$ is representable exactly as a mixture of the above mentioned beta densities. The mixture model is especially useful for Bayesian inference, where priors can be induced upon the mixture densities through a Dirichlet process prior on the mixing distribution. When hypotheses are randomly assigned as null or alternative with a specific probability, then the p-value distribution is a mixture of a uniform component and a mixture of beta densities of the type mentioned above. By applying the general theory of posterior consistency for density estimation, we have shown that the posterior distribution for estimating the density of p-values is consistent at the true density if it is of the given form and the prior on the mixing distribution has every distribution in its weak support. Under some further conditions which essentially separate mixtures of beta densities from the uniform, it follows that posterior consistency for density estimation leads to consistency in estimating positive false discovery rates for multiple hypotheses testing. This property gives asymptotic justification of a recently proposed Bayesian method of estimating positive false discovery rates by the same set of authors.